\documentclass{amsart}
\usepackage{amscd,amssymb,graphics,amsthm,amsmath,amssymb}

\theoremstyle{plain}
\newtheorem{theorem}{Theorem}

\newtheorem{lemma}[theorem]{Lemma}

\newtheorem*{main}{Main Theorem}

\theoremstyle{definition}
\newtheorem{definition}{Definition}

\theoremstyle{remark}

\numberwithin{equation}{section}

\begin{document}

\title[The Morse Complex on a Manifold With Corners]{The Morse Complex for a 
Morse Function on a Manifold with Corners}
\author{David G.C. Handron}
\date{\today}
\address{Department of Mathematical Sciences\\
Carnegie Mellon University \\
	Pittsburgh, PA \\
	(412)268-2545}
\email{handron@andrew.cmu.edu}
\date{\today}

\begin{abstract}
A Morse function $f$ on a manifold with corners $M$ allows the
characterization of the Morse data for a critical point by the Morse
index. In fact, a modified gradient flow allows a proof of the Morse
theorems in a manner similar to that of classical Morse theory. It
follows that $M$ is homotopy equivalent to a CW-complex with one cell
of dimension $\lambda$ for each essential critical point of index
$\lambda$. The goal of this article is to determine the boundary maps
of this CW-complex, in the case where $M$ is compact and orientable.

First, the boundary maps are defined in terms of the modified gradient
flow.  Then a transversality condition is imposed which insures that
the attaching map is non-degenerate in a neighborhood of each critical
point. The degree of this map is then interpreted as a sum of
trajectories connecting two critical points each counted with a
multiplicity determined by a choice of orientations on the tangent
spaces of the unstable manifold at each critical point.
\end{abstract}

\maketitle

\section{Introduction}
Goresky and MacPherson's {\em Stratified Morse Theory}
(\cite{Gor-Mac}) represented a great step forward in extending the
ideas of Morse, Thom and Smale to spaces more general than smooth manifolds.
The complexities involved in dealing with stratified spaces prompted
Goresky and MacPherson to comment on the nostalgia their Stratified
Morse Theory might inspire for the classical version, where the Morse
data of a critical point is determined by a single number, the Morse
index.

In \cite{Vakhrameev} Vakhremeev proved the Morse theorems in the
setting of manifolds with corners, a setting which strikes a nice
balance between the simplicity of classical Morse theory on one hand,
and the generality of Stratified Morse Theory on the other.  Manifolds
with corners are a class of stratified space that arise naturally in
many applications, yet on these spaces the Morse data for an
(essential) critical point still is determined by a single number.

Others have studied Morse Theory in similar settings, (e.g. Hamm
(\cite{Hamm2},
\cite{Hamm2.5}), Siersma (\cite{Siersma}) and Braess (\cite{Braess}), 
but none of these make use of a global flow on manifolds with corners.
The Morse theory for manifolds with corners developed in
\cite{handron} utilizes a modified gradient flow to prove the theorems
in a more classical manner (e.g. as in~\cite{Milnor}).  This allows,
as we shall see, the construction of a Morse complex with boundary
maps determined by the trajectories connecting (essential) critical
points.

\subsection{Setup and Definitions}

Let $\{ {\mathbf e}_1,\dots,{\mathbf e}_\ell \}$ denote the standard
basis vectors in ${\mathbb R}^n$.  Define ${\mathbb H}_\ell^n$ to be the
set
$$
{\mathbb H}_\ell^n=\{ {\mathbf w}\in{\mathbb R}^n: {\mathbf
w}\cdot{\mathbf e}_i\geq 0 \text{ for all }i>n-\ell \},
$$
where $\cdot$ denotes the standard inner product on ${\mathbb R}^n$.

\begin{definition} An {\em $n$-dimensional manifold with corners} is a
topological space $M$ together with an atlas ${\mathcal{A}}$ of charts
${\mathbf x}_a:U_a\to{\mathbb H}^n_{\ell_a}$ such that 
$\{U_a\}_a$ is an open cover of $M$.  We say that $(M,\mathcal{A})$ is
a $C^r$-manifold with corners if each transition function
$\mathbf{x}_b^{}\mathbf{x}_a^{-1}:\mathbf{x}_a(U_a\cap
U_b)\to\mathbf{x}_b(U_a\cap U_b)$ is of class $C^r$, i.e. can be
extended to a $C^r$ map between open sets of $\mathbb{R}^n$. 
\label{mwc}
\end{definition}

The tangent space $T_pM$ of a manifold with corners $M$ at a point $p$ may be
thought of as equivalence classes of the set
$$
\{ ({\mathbf x},{\mathbf v}):{\mathbf x}\text{ is a coordinate chart at }
p\in M\text{ and }{\mathbf v}\in{\mathbb R}^n \},
$$
where $({\mathbf x},{\mathbf v})\sim({\mathbf y},{\mathbf w})$ if
$D({\mathbf x}\circ{\mathbf y}^{-1})({\mathbf w})={\mathbf v}$.  
Then for each $p\in M$, $T_pM$ is an $n$-dimensional vector space.

If $p\in M$, we will say a coordinate chart {\em at $p$} is a chart
${\mathbf x}_p\in{\mathcal{A}}$ such that ${\mathbf
x}_p(p)=\mathbf{0}\in{\mathbb H}^n_\ell$.  In this case, the number
$\ell=\ell(p)$ is uniquely determined by $p$.  This allows
consideration of $M$ as a stratified space.  Each connected component
of ${\mathcal E}_\ell=\{p\in M:\ell(p)=\ell\}$ is a stratum of
dimension $n-\ell$.

\begin{definition}
A tangent vector in $T_pM$ points {\em outward} if some representative
$({\bold x},{\bold v})$ has ${\bold v}\notin{\mathbb H}^n_{\ell(p)}$.  A tangent 
vector is $T_pM$ points {\em inward} (or {\em into} M) if some representative 
$({\bold x},{\bold v})$ has ${\bold v}\in{\mathbb H}^n_{\ell(p)}$.
\label{in_out_vect}
\end{definition}

These terms are well defined, since for any two coordinate charts at
$p$, the transition functions preserve ${\mathbb H}^n_{\ell(p)}$.
Note that the definition of an inward pointing vector includes those
vectors which are tangent to the boundary of $M$.

If $K$ is a stratum of $M$ then $K$ is a manifold.  We may thus
consider at a point $p\in K$ the tangent space $T_pK$.  In fact,
$T_pK\subseteq T_pM$ as
$$
T_pK=\{\mathbf{v}\in T_pM|\mathbf{v}\mbox{ and }-\mathbf{v}\mbox{ point inward}\}
$$

If $p\in K$ is in the closure of another stratum, $L$, we can define the 
generalized tangent space 
$$
T_pL=\{ {\bold w}\in T_pM: {\bold w}=\lim_{i\to\infty}{\bold v}_i\in T_{q_i}K
\text{ for some sequence }\{q_i\}\to p \}.
$$
We may also write this as $T_pL=\lim_{q\to p}T_qL$.

\begin{definition}
A {\em Riemannian metric} on a $C^r$-manifold with corners $M$
consists of an inner product $\langle\cdot,\cdot\rangle_p$ on each
tangent space $T_pM$, such that with respect to any coordinate chart
$\mathbf{x}:U\to\mathbb{H}^n_\ell$, the functions
$\langle\frac{\partial}{\partial x_i},\frac{\partial}{\partial
x_j}\rangle:\mathbf{x}(U)\to\mathbb{R}$ are of class $C^r$, i.e. can
be extended to a $C^r$ function on an open set in $\mathbb{R}^n$.
\end{definition}

We say that a coordinate chart ${\mathbf x}:U\to \mathbb{H}^n_\ell$ at a
point $p$ is a {\em standard} coordinate chart if whenever $i\leq n-\ell$
and $j>n-\ell$ the tangent vectors $\frac{\partial}{\partial x_i}$ and
$\frac{\partial}{\partial x_j}$ are orthogonal.

\begin{lemma}
For each point $p\in M$, there is a standard coordinate chart at $p$.
\end{lemma}

\begin{proof}
Let
$\mathbf{y}:U\to\mathbb{H}^n_\ell\cong\mathbb{R}^{n-\ell}\times[0,\infty)^\ell$
be a coordinate chart at $p$.  Then $\mathbf{y}(p)=\mathbf{0}$.
Define an inner product on $\mathbb{R}^n$ by
$$
\langle\mathbf{v}, \mathbf{w}\rangle =
\langle\mathbf{y}_\ast^{-1}(\mathbf{v}),
\mathbf{y}_\ast^{-1}(\mathbf{w})\rangle_p.
$$
Here we have used the natural identification of $\mathbb{R}^n$ with
$T_\mathbf{0}\mathbb{R}^n$.

Using this inner product, define an orthogonal projection onto
$\mathbb{R}^{n-\ell}$.  For each vector $\mathbf{v}\in\mathbb{R}^n$,
$\mbox{proj}(\mathbf{v})\in\mbox{span}(\mathbf{e}_1,\dots,\mathbf{e}_{n-\ell})$
is the unique vector such that
$\langle\mathbf{v}-\mbox{proj}(\mathbf{v}),\mathbf{w}\rangle=0$ for
every
$\mathbf{w}\in\mbox{span}(\mathbf{e}_1,\dots,\mathbf{e}_{n-\ell})$.

Now define a new basis for $\mathbb{R}^n$ by
$$
\begin{array}{rlc}
\overline{\mathbf{e}}_i\negmedspace\negthickspace\negthickspace &=\mathbf{e}_i & \mbox{for }i\leq n-\ell \\
\overline{\mathbf{e}}_j\negmedspace\negthickspace\negthickspace &=\mathbf{e}_j-\mbox{proj}(\mathbf{e}_j) & j>n-\ell.
\end{array}
$$
Then if $i\leq n-\ell$ and $j>n-\ell$,
$\langle\overline{\mathbf{e}}_i,\overline{\mathbf{e}}_j\rangle=0$.

Let $J:\mathbb{R}^n\to\mathbb{R}^n$ be the linear transformation that
carries each $\mathbf{e}_k$ to $\overline{\mathbf{e}}_k$.  Then $J$ is
an isomorphism.  We wish to show that
$J:\mathbb{H}^n_\ell\to\mathbb{H}^n_\ell$.  It is sufficient to show
that if
$$
\mathbf{v}=a_1\mathbf{e}_1+\dots+a_n\mathbf{e}_n
$$
with $a_j\geq0$ whenever $j>n-\ell$, then the expression
$$
\mathbf{v}=\overline{a}_1\overline{\mathbf{e}}_1+\dots+\overline{a}_n\overline{\mathbf{e}}_n
$$
satisfies $\overline{a}_j\geq0$ whenever $j>n-\ell$.  But if
$\mathbf{v}=\overline{a}_1\overline{\mathbf{e}}_1+\dots+\overline{a}_n\overline{\mathbf{e}}_n$, then
$$
\begin{array}{rl}
\mathbf{v}&=\overline{a}_1\mathbf{e}_1+\dots+\overline{a}_{n-\ell}\mathbf{e}_{n-\ell} \\
&\quad\quad+\overline{a}_{n-\ell+1}\left(\mathbf{e}_{n-\ell+1}-\mbox{proj}(\mathbf{e}_{n-\ell+1})\right)+\dots+
\overline{a}_{n}\left(\mathbf{e}_{n}-\mbox{proj}(\mathbf{e}_{n})\right). \\
\end{array}
$$
Then, setting
$\mbox{proj}(\mathbf{e}_j)=\sum_{k=1}^{n-\ell}p_{j,k}\mathbf{e}_k$, we
see that
$$
\begin{array}{rl}
\mathbf{v}&=\left(\overline{a}_1-\sum_{j=n-\ell+1}^{n}p_{j,1}\right)\mathbf{e}_{1}+\dots+\left(\overline{a}_{n-\ell}-\sum_{j=n-\ell+1}^{n}p_{j,n-\ell}\right)\mathbf{e}_{n-\ell}
\\
&\quad\quad +\overline{a}_{n-\ell+1}\mathbf{e}_{n-\ell+1}+\dots+\overline{a}_{n-\ell}\mathbf{e}_{n-\ell}.
\end{array}
$$
So for $j>n-\ell$, $\overline{a}_{j}=a_j$.  It follows that
$J:\mathbb{H}^n_\ell\to\mathbb{H}^n_\ell$.

Finally, let $\mathbf{x}=J\circ\mathbf{y}$.  Then
$\mathbf{x}:U\to\mathbb{H}^n_\ell$ and $\mathbf{x}(p)=\mathbf{0}$, so
$\mathbf{x}$ is a standard coordinate chart at $p$.  Moreover,
$$
\mathbf{x}_\ast^{-1}\left(\frac{\partial}{\partial
x_k}\right)=J\circ\mathbf{y}_\ast^{-1}\left(\frac{\partial}{\partial
y_k}\right).
$$  
Since
$$
\mathbf{y}_\ast^{-1}\left(\frac{\partial}{\partial
y_k}\right)=\mathbf{e}_k,
$$
 it follows that
$\mathbf{x}_\ast^{-1}\left(\frac{\partial}{\partial
x_k}\right)=\overline{\mathbf{e}}_k$, and $\mathbf{x}$ is a standard
coordinate chart at $p$.

\end{proof}

Let $M$ be a compact, $n$-dimensional $C^r$-manifold with corners
($r\geq 2$) endowed with a Riemannian metric.  Let $f:M\to\mathbb{R}$
be a $C^r$ function.  For every stratum $K$ of $M$, the function
$f|_K:K\to\mathbb{R}$ is a $C^r$ function on the manifold $K$.  If
$p\in K$, we say $p$ is a {\it critical point} of $f$ whenever $p$ is
a critical point of $f\vert_K$. Recall that a critical point of $f|_K$
is {\em non-degenerate} if the Hessian of $f|_K:K\to\mathbb{R}$ is
non-zero at $p$.

If $p$ is an interior point of $M$, then it is contained in a stratum
$K$ with dimension $n$.  The gradient $\nabla f$ of the function
$f:M\to\mathbb{R}$ has a natural definition.  $K$ is a manifold with a
Riemannian metric induced by the metric on $M$.  We let $\nabla
f(p)=\nabla f|_K(p)\in T_pK=T_pM$.  If $p$ is in the boundary of $M$, we
can define $\nabla f(p)=\lim_{q\to p}\nabla f(q)$ where $q$ is taken
to be in the interior of $M$.

We can now define the class of functions we wish to consider:

\begin{definition}
Let $f:M\to{\mathbb R}$ be a $C^r$ function from a compact,
$n$-dimensional $C^r$-manifold with corners ($r\geq 2$) endowed with a
Riemannian metric into the real numbers.  Then $f$ is a Morse function
on if the following hold:

\begin{enumerate}
	\item If $K$ is a stratum of $M$, and $p\in K$ is a critical point of 
		$f\vert_K:K\to {\mathbb R}$, then either 
\begin{enumerate}
		\item\label{onea} $p$ is a non-degenerate critical point of 
			$f\vert_K:K\to {\mathbb R}$ or
		\item\label{oneb} the vector $-\nabla f(p)$ points into $M$.
\end{enumerate}
	\item If $p\in K$ is a critical point of $f|_K$, then for any stratum
		$L$ with $p$ in the closure of $L$, $-\nabla
		f(p)\notin T_pL^\perp\subseteq T_pM$.  

	\item\label{three} If $-\nabla f$ is tangent to a stratum $K$
		at a point $p\in K$, and
		$\mathbf{x}:U\to\mathbb{H}^n_\ell$ is a standard
		coordinate chart at $p$, then for $j>n-\ell$, the
		directional derivative of $\frac{\partial f}{\partial
		x_j}$ in the direction $-\nabla f(p)$ is non-zero.
\end{enumerate}
\end{definition}

This definition is inspired by the definition of Morse function in
\cite{Gor-Mac}.  It differs in condition~(\ref{oneb}), where the
critical point is not required to be non-degenerate, and in the
addition of condition~(\ref{three}), which ensures that the vector
field $G$ defined below produces a continuous flow.

We say that $p$ is an {\em essential critical point} of $f$ if $f$ satisfies 
condition~(\ref{onea}) at $p$, but not condition~(\ref{oneb}).  The {\em index} of 
an essential critical point $p$ in a stratum $K$ is equal to its index (in
the classical sense) as 
a critical point of $f|_K$.

Suppose $f:M\to{\mathbb R}$ is a Morse function such that 
$f^{-1}(-\infty,r]$ is compact for each $r\in\mathbb{R}$.  It is shown in 
\cite{handron} that $M$ is homotopy equivalent to a CW-complex with one 
cell of dimension
$\lambda$ for each essential critical point of
index $\lambda$. 

Let $G$ be the {\em modified gradient vector field} of $f$ on $M$,
defined in \cite{handron} by projecting the vector $-\nabla f(p)$ onto
the stratum of highest dimension such that the resulting vector does
not point outward from $M$.  It is shown in~\cite{handron} that $G$
induces a flow $\varphi:[0,\infty)\times M\to M$ and that the
stationary points of this flow are exactly the essential critical
points of $f$.  The material above can be found in greatly expanded
form in \cite{handron}.

We can define the stable and 
unstable sets of an essential critical point $p$ by
$$
S(p)=\left\{q\in M:\lim_{t\to\infty}\varphi(t,q)=p\right\},
$$
and
$$
U(p)=\left\{q\in M:\exists \{q_k\}_{k=1}^\infty\subset M \mbox{ such that } 
\varphi(k,q_k)=q \mbox{ and }\lim_{k\to\infty}q_k=p\right\}.
$$

\section{Results Proved}

Our goal is to describe the boundary maps of the CW-complex generated
by $f$ in terms of the trajectories of the flow $\varphi$.  If we
assume a Smale-like condition (Definition~\ref{smale_condition}) for
the function $f$ then the number of trajectories connecting a pair of
critical points whose indeces differ by one is finite.  Furthermore,
if $M$ is orientable and the index of $p_2$ is one greater than that
of $p_1$, then we can associate to each trajectory $\tau$ connecting
$p_2$ to $p_1$ a sign $+1$ or $-1$.  This sign is determined by
comparing, via the flow $\varphi$, orientations on $U(p_1)$ and
$U(p_2)$.  The degree of the pair $(p_2,p_1)$ is
$$
\mbox{degree}(p_2,p_1)=\sum_{\tau}\mbox{sign}(\tau),
$$
where the sum is over all trajectories connecting $p_2$ to $p_1$.

If $P_1$ and $P_2$ are the cells in the CW-complex corresponding to the 
critical points $p_1$ and $p_2$, then the degree of the attaching map of $P_2$
along $P_1$ is the sum of the signs of the trajectories from $p_2$ to $p_1$.

Now suppose $M$ is a compact oriented manifold with corners and
$f:M\to{\mathbb R}$ is a Morse function satisfying the conditions of
Definition~\ref{smale_condition}.  The goal of this article is to
define a chain complex which is determined by the function $f$ and
whose homology is that of $M$.  Let $C_\ell(f)$ be the set of essential
critical points of $f$ with index $j$, and $V_\ell$ the free Abelian
group of formal ${\mathbb Z}$-linear combinations of the elements of
$C_\ell(f)$.  We can define $\partial_\ell:V_\ell\to V_{j-1}$ by setting
$$
\partial_\ell(p)=\sum_{q\in C_{j-1}(f)}\mbox{degree}(p,q)q
$$
for each $p\in C_\ell(f)$, and extending linearly to all of $V_\ell$.
Finally, we can state the

\begin{main}
\textbf{(Theorem~\ref{main_theorem})}
The free Abelian groups $V_\ell$ and maps $\partial_\ell$ form a chain complex whose 
homology groups are identical to the $\mathbb{Z}$-homology groups 
 of the topological space $M$.
\end{main}

\section{The Attaching Map of a Cell in the Morse Complex}

\label{attaching}

Let $p_1$ and $p_2$ be essential critical points with indices
$\lambda_1$ and $\lambda_2$ respectively.  Let $c_i=f(p_i)$ for
$i=1,2$.  We can choose $\varepsilon$
sufficiently small that the set 
$$
U_\varepsilon(p_i)=U(p_i)-f^{-1}((c_i-\varepsilon])=U(p_i)-M_{c_i-\varepsilon}
$$ 
is homeomorphic to the interior of $e^{\lambda_i}$, a cell with
dimension $\lambda_i$, for $i=1,2$.  Here, we follow the convention
that $M_{a}=f^{-1}((-\infty,a])$.  In addition,
suppose that $f(p_1)=c_1$ and $f(p_2)=c_2$, where $c_1<c_2$, and that
there are no other critical points in $f^{-1}([c_1,c_2])$.

The
flow $\varphi$ can be viewed as a deformation retraction, and used to
define homotopies
$$
M_{c_2-\varepsilon}\simeq\varphi(t,M_{c_2-\varepsilon})
$$
and
\begin{equation}
\label{first_case}
M_{c_2-\varepsilon}\cup U(p_2)\simeq\varphi(t,M_{c_2-\varepsilon}\cup U(p_2))
\end{equation}

Since $U(p_2)$ is invariant under the flow $\varphi$ , it follows that
\begin{equation}
\varphi(t,M_{c_2-\varepsilon}\cup
U(p_2))=\varphi(t,M_{c_2-\varepsilon})\cup U(p_2).
\label{middle_case}
\end{equation}
The space on the left hand side of (\ref{first_case}) is
$M_{c_2-\varepsilon}$ with a $\lambda_2$-cell attached, but on the
right hand side of (\ref{middle_case}),
$U(p_2)-\varphi(t,M_{c_2-\varepsilon})$ may not be homoemorphic to an
open cell.  This is a notable difference from the classical theory.

To simplify our thinking, we can consider
$M_{c_2-\varepsilon}\cup U(p_2)$ to be $M_{c_2-\varepsilon}$ with a
$\lambda_2$-cell attached via the inclusion map
$$
\sigma:\partial U_\varepsilon(p_2)\hookrightarrow M_{c_2-\varepsilon}.
$$
Then the family of homotopy equivalences
$$
M_{c_2-\varepsilon}\simeq \varphi(t,M_{c_2-\varepsilon}),
$$
gives rise to homotopy equivalences
\begin{equation}
\label{second_case}
M_{c_2-\varepsilon}\cup_\sigma e^{\lambda_2} \simeq
\varphi(t,M_{c_2-\varepsilon})
\cup_{\varphi(t,\cdot)\circ\sigma}e^{\lambda_2}.
\end{equation}
Here we have identified $\overline{U_\varepsilon(p_2)}$ with
$e^{\lambda_2}$.

In the first case (\ref{middle_case}), we allow the whole space to
flow by $\varphi$.  The equality follows from the $\varphi$-invariance
of $U(p_2)$.  In the second case (\ref{second_case}), we simply allow
$M_{c_2-\varepsilon}$ to flow by $\varphi$, and with it the image of
the attaching map $\sigma_2$.  (c.f. Lemma~3.7 in \cite{Milnor}.)

From here it is a short step to prove 

\begin{lemma}
Let $p_1$ and $p_2$ be two essential critical points of $f$ such that 
$f(p_1)=c_1$, $f(p_2)=c_2$, $c_1<c_2$ and there are no other critical 
points in $f^{-1}([c_1-\varepsilon,c_2])$.  Then there is a map
$\Phi_2:M_{c_2-\varepsilon}\to M_{c_1+\varepsilon}$ such that
$$
M_{c_2-\varepsilon}\cup U(p_2)=M_{c_2-\varepsilon}\cup_\sigma
e^{\lambda_2}\simeq M_{c_1+\varepsilon}\cup_{\Phi_2\circ\sigma} e^{\lambda_2},
$$
where $\lambda_2$ is the index of $p_2$.  Moreover, $\Phi_2$ may be
defined in a way that preserves $\varphi$-trajectories in the sense
that for every $q\in M_{c_2-\varepsilon}$, the points $q$ and
$\Phi_2(q)$ lie on the same $\varphi$-trajectory.
\label{lemma_1}
\end{lemma}

\begin{proof}
The equality follows immediately from the discussion above.  To define
the map $\Phi_2$, we use the flow $\varphi$ to construct a homotopy.
For each $q\in M_{c_2-\varepsilon}$, define 
$$
t_q^{(1)}=\inf\left\{t\in\mathbb{R}^+:\varphi(t,q)\in M_{c_1+\varepsilon}\right\}
$$
The function $q\mapsto t_q^{(1)}$ is continuous. Since there are no
critical values between $c_1$ and $c_2$, $t_q^{(1)}$ is finite for each 
$q\in M_{c_2-\varepsilon}$.  Moreover, since $M_{c_2-\varepsilon}$ is
compact,
$T_1=\sup\{t_q^{(1)}:q\in M_{c_2-\varepsilon}\}$ is finite.

Now we can define a  homotopy 
$H_1:[0,1]\times M_{c_2-\varepsilon}\to M_{c_1-\varepsilon}$
by
$$
H_1(s,q)=\left\{
\begin{array}{cc}
\varphi(T_1s,q)	& T_1s\leq t_q^{(1)} 	\\
\varphi(t_q^{(1)},q)	& T_1s\geq t_q^{(1)}
\end{array}
\right.
$$

If we set $\Phi_2(q)=H_1(1,q)$ for $q\in M_{c_2-\varepsilon}$, the
homotopy equivalence in the lemma follows.
\end{proof}

$\Phi_2$ can be interpreted as allowing each point to flow by
$\varphi$ until it reaches $M_{c_1+\varepsilon}$. 
Now we must deal with the more complicated situation in the
neighborhood of the critical point $p_1$.

\begin{lemma}
Let $p_1$ be an essential critical point of $f$ such that $f(p_1)=c_1$
and let $\varepsilon$ be small enough that there are no other critical
points in $f^{-1}([c_1-\varepsilon,c_1+\varepsilon])$.  Then there is
a homotopy equivalence
$$
\Psi_2:M_{c_1+\varepsilon}\to M_{c_1-\varepsilon}\cup U(p_1).
$$  
Moreover, $\Psi_2$ can be chosen such that $\Psi_2(q)=p_1$ if and only
if $q\in S(p_1)\cap M_{c_1+\varepsilon}$.
\label{lemma_2}
\end{lemma}

\begin{proof}
Since $p_1$ is an essential critical point with index $\lambda_1$, we
can choose a coordinate chart such that
$$
f(\mathbf{x})=c_1-x_1^2-\dots-x_{\lambda_1}^2+x_{\lambda_1+1}^2+\dots+x_{n-j}^2+
x_{n-j+1}+\dots+x_n
$$
We can define, for $\varepsilon$ sufficiently small, a region $H$ and
show that $M_{c_1-\varepsilon}\cup H$ is homotopy equivalent to both
$M_{c_1+\varepsilon}$ and $M_{c_1-\varepsilon}\cup U(p_1)$.  To do so
we first choose a smooth function $\mu:\mathbb{R}\to\mathbb{R}$ so
that $\mu(0)>\varepsilon$, $\mu(r)=0$ for $r>2\varepsilon$, and
$-1<\mu'\leq0$.  Then we set
$$
F(\mathbf{x})=f(\mathbf{x})-\mu(x_1^2+\dots+x_{\lambda_1}^2+2(x_{\lambda_1+1}^2+\dots+x_{n-j}^2)+2(x_{n-j+1}+\dots+x_n))
$$
Finally, we define
$H=\overline{F^{-1}(-\infty,c_1-\varepsilon]-M_{c_1-\varepsilon}}$.
Note that $H$ is used to denote this space, as it is in \cite{Milnor},
despite the possibility of confusion with the homotopies defined above
and below.  In this article, $H_i$ is a homotopy and $H$ is a subset
of $M$.

The proof of the lemma will be in two parts.  First, we construct a
homotopy to show that $M_{c_1+\varepsilon}\simeq
M_{c_1-\varepsilon}\cup H$.  We then construct a different homotopy
to show that $M_{c_1-\varepsilon}\cup H\simeq M_{c_1-\varepsilon}\cup
U(p_1)$.

For each $q\in M_{c_1+\varepsilon}$, define 
$$
t_q^{(2)}=\inf\left\{t\in\mathbb{R}^+:\varphi(t,q)\in M_{c_1-\varepsilon}\cup H\right\}
$$
The function $q\mapsto t_q^{(2)}$ is continuous, and $t_q^{(2)}$ is finite for each 
$q\in M_{c_2-\varepsilon}$.  Moreover, since $M_{c_2-\varepsilon}$ is compact,
$T_2=\sup\{t_q^{(2)}:q\in M_{c_2-\varepsilon}\}$ is finite.

Now we can define the first of our two homotopies 
$H_2:[0,1]\times M_{c_1+\varepsilon}\to M_{c_1-\varepsilon}\cup H$
by
$$
H_2(s,q)=\left\{
\begin{array}{cc}
\varphi(T_2s,q)	& T_2s\leq t_q^{(2)} 	\\
\varphi(t_q^{(2)},q)	& T_2s\geq t_q^{(2)}
\end{array}
\right.
$$
This shows that  $M_{c_1+\varepsilon}= H_2(0,M_{c_1+\varepsilon})\simeq
H_2(1,M_{c_1+\varepsilon})= M_{c_1-\varepsilon}\cup H$.

Now we show that $M_{c_1-\varepsilon}\cup H\simeq
M_{c_1-\varepsilon}\cup U(p_1)$.  In our chosen coordinate system the
stable set of $p_1$ corresponds to
$$
\{(x_1,\dots,x_n)|x_1=\dots=x_\lambda=0\}.
$$
  The unstable set of
$p_1$ corresponds to 
$$
\{(x_1,\dots,x_n)|x_{\lambda+1}=\dots=x_n=0\}.
$$

Consider the function 
%$g(\mathbf{x})=f(p_1)+x_{\lambda+1}^2+\dots+x_{n-j}^2+x_{n-j+1}$.
$g(\mathbf{x})=f(p_1)+x_{\lambda+1}^2+\dots+x_{n-j}^2+x_{n-j+1}+\dots+x_n$.
This function is constant on the unstable set, with a value of
$f(p_1)$.  For $q$ outside the unstable set, $f(q)>f(p_1)$.  Moreover,
on the stable set the $-\nabla g$ agrees with $-\nabla f$.

Let $X$ denote the modified (negative) gradient vector field of the function
$g$.  Then $X(q)=G(q)$ for $q\in S(p_1)$, and $X(q)=0$ for $q\in U(p_1)$.
The flow of $X$ is directed toward $U(p_1)$.  Let $\psi$ denote the flow of 
the vector field $\frac{X}{\|X\|}$, and note that $\psi$ preserves the 
$\varphi$-trajectories that have $p_1$ as a lower endpoint.

For each $q$ in $H$, we can set 
$t_q^{(3)}=\inf\{t\in\mathbb{R}^+:\psi(t,q)\in M_{c_1-\varepsilon}\cup U(p_1)\}$, and for $q\in M_{c_1-\varepsilon}$ set $t_q^{(3)}=0$.
Then define $T_3=\sup\{t_q^{(3)}:q\in H\}$.
Then the homotopy 
$H_3:[0,1]\times M_{c_1-\varepsilon}\cup H\to 
M_{c_1-\varepsilon}\cup U(p_1)$  
given by
$$
H_3(s,q)=\left\{
\begin{array}{cc}
\psi(T_3s,q)	& T_3s\leq t_q^{(3)}	\\
\psi(t_q^{(3)},q)	& T_3s\geq t_q^{(3)}
\end{array}
\right.
$$
gives a homotopy equivalence 
$H_3(1,\cdot):M_{c_1-\varepsilon}\cup H\to M_{c_1-\varepsilon}\cup U(p_1)$.

Then the composition $\Psi_2=H_3(1,\cdot)\circ H_2(1,\cdot)$ is a
homotopy equivalence $\Psi_2:M_{c_1+\varepsilon}\to
M_{c_1-\varepsilon}\cup U(p_1)$.

\end{proof}

There is no reason, in principle, that the $\varepsilon$ used to
define $U_\varepsilon(p_2)$ must be the same as the $\varepsilon$ in
Lemma~\ref{lemma_1} and Lemma~\ref{lemma_2}.  Nor must the
$\varepsilon$'s from the two Lemmas agree with each other.  They can
be chosen to be the same, however, and this serves to simplify the
notation in the remainder of our argument.

Now lets view this in terms of attaching maps.  Let $\sigma_2$ denote
the inclusion map $\sigma_2:\partial U_\varepsilon(p_2)\to
M_{c_2-\varepsilon}$.  (Previously this was simply $\sigma$.)  Then we
can use the homotopy equivalence $\Phi_2$ to define an attaching map
$$
\Phi_2\circ\sigma_2:\partial U_\varepsilon(p_2)\to M_{c_1+\varepsilon}.
$$
Now, identifying $U_\varepsilon(p_2)$ with $e^{\lambda_2}$, we find that
$$
M_{c_2+\varepsilon}\simeq M_{c_2-\varepsilon}\cup U(p_2)\simeq
M_{c_1+\varepsilon}\cup_{\Phi_2\circ\sigma_2}e^{\lambda_2}.
$$
The first homotopy equivalence follows from Lemma~\ref{lemma_2}.

Using the homotopy equivalence $\Psi_2$, we can define an attaching map
$$
\Psi_2\circ\Phi_2\circ\sigma_2:\partial U_\varepsilon(p_2)\to
M_{c_1-\varepsilon}\cup U(p_1).
$$
From this, it follows that
$$
M_{c_2+\varepsilon}\simeq M_{c_2-\varepsilon}\cup U(p_2)\simeq
M_{c_1-\varepsilon}\cup
U(p_1)\cup_{\Psi_2\circ\Phi_2\circ\sigma_2}e^{\lambda_2}.
$$

Let $\Sigma_2=\Psi_2\circ\Phi_2\circ\sigma_2$.  Because $\Phi_2$
preserves all the $\varphi$-trajectories and $\Psi_2$ preserves those
in $S(p_1)$, it follows that $\Sigma_2(q)=p_1$ if and only if $q\in
\partial U_\varepsilon(p_2)$ lies on a trajectory connecting $p_2$ to $p_1$.

Using these ideas, we can construct a cellular complex (which will
later prove to be a CW-complex) that is homotopy equivalent to $M$.
Let $p_1,p_2,\dots,p_m$ be the essential critical points of $M$, each
$p_i$ having index $\lambda_i$.  Assume the values $f(p_i)$ are
distinct and the essential critical points are labeled such that
$f(p_i)<f(p_j)$ whenever $i<j$.  Then $c_1=f(p_1)$ must be an absolute
minimum, and $M_{c_1+\varepsilon}$ is homotopy equivalent to the
0-cell, $U_\varepsilon(p_1)$.  In fact the map $\Psi_2$ is a homotopy
equivalence.

Now, as in Lemma~\ref{lemma_1} and Lemma~\ref{lemma_2}, we consider the 
inclusion map 
$\sigma_2:\partial U_\varepsilon(p_2)\hookrightarrow M_{c_2-\varepsilon}$.
Viewing $U_\varepsilon(p_2)$ as a $\lambda_2$-cell,
and using the attaching map 
$\Sigma_2=\Psi_2\circ\Phi_2\circ\sigma_2$, of
$U_\varepsilon(p_2)$ onto $U_\varepsilon(p_1)=\{p_1\}$ 
we see that 
$$
\begin{array}{rl}
M_{c_2+\varepsilon} &\simeq
U_\varepsilon(p_1)\cup_{\Sigma_2}U_\varepsilon(p_2) \\
& \simeq e^{\lambda_1}\cup_{\Sigma_2}e^{\lambda_2}.
\end{array}
$$

If $\varepsilon$ is sufficiently small, then $U_\varepsilon(p_3)$ is
homeomorphic to a $\lambda_3$-cell.  It takes only a slight extension
to define the attaching map for this cell.  Let $\sigma_3$ be the
inclusion map $\partial U_\varepsilon(p_3)\hookrightarrow
M_{c_3-\varepsilon}$, and define $\Phi_3$, $\Psi_3$ for the critical
points $p_2$ and $p_3$ using Lemma~\ref{lemma_1} and
Lemma~\ref{lemma_2}.  Then
$\widetilde{\Sigma}_3=\Psi_3\circ\Phi_3\circ\sigma_3$ gives the
appropriate attaching map of $\partial U_\varepsilon(p_3)$ onto
$M_{c_2-\varepsilon}\cup U(p_2)\simeq
M_{c_2-\varepsilon}\cup_{\sigma_2}e^{\lambda_2}$.

Since $\Psi_2\circ\Phi_2$ is a homotopy equivalence from
$M_{c_2-\varepsilon}\cup_{\sigma_2}e^{\lambda_2}$ to
$e^{\lambda_1}\cup_{\Sigma_2}e^{\lambda_2}$, the desired attaching map
is given by
$$
\Sigma_3=\Psi_2\circ\Phi_2\circ\Psi_3\circ\Phi_3\circ\sigma_3.
$$
It follows that
$$
\begin{array}{rl}
M_{c_3+\varepsilon} &\simeq
U_\varepsilon(p_1)\cup_{\Sigma_2}U_\varepsilon(p_2)\cup_{\Sigma_3}U_\varepsilon(p_3) \\
& \simeq e^{\lambda_1}\cup_{\Sigma_2}e^{\lambda_2}\cup_{\Sigma_3}e^{\lambda_3}.
\end{array}
$$

Because of the construction of $\Phi_3$ and $\Psi_3$,
$\Sigma_3(q)=p_2$ if and only if $q\in\partial U_\varepsilon(p_3)$
lies on a $\varphi$-trajectory connecting $p_3$ to $p_2$.  We would
like to be able to say that $\Sigma_3(q)=p_1$ if and only if
$q\in\partial U_\varepsilon(p_3)$ lies on a $\varphi$-trajectory
connecting $p_3$ to $p_1$.  This will be true provided that $\Psi_3$
does not disrupt these trajectories.  

This can not be true in every case.  With care, though, we can be sure
that trajectories connecting critical points whose indices differ by
one are not disturbed.  This issue is resolved by
Lemmas~\ref{critical_point_lemma} and \ref{cw-complex_lemma}.

In general, for each critical point $p_k$ we define
$$
\Sigma_k=\Psi_2\circ\Phi_2\circ\dots\Psi_k\circ\Phi_k\circ\sigma_k.
$$
These maps gives a homotopy equivalence
$$
\begin{array}{rl}
M & \simeq
U_\varepsilon(p_1)\cup_{\Sigma_2}U_\varepsilon(p_2)\cup\dots\cup_{\Sigma_m}U_\varepsilon(p_m)
\\
& \simeq e^{\lambda_1}\cup_{\Sigma_2}e^{\lambda_2}\cup\dots\cup_{\Sigma_m}e^{\lambda_m}.
\end{array}
$$

We shall see that under appropriate conditions 
(Definition~\ref{smale_condition}), this cellular complex is in fact
a CW-complex.

\section{The Degree of an Attaching Map}
\label{degree}

Our goal is to show that, as in classical Morse theory, the degree of
this attaching map can be determined by counting the trajectories
connecting essential critical points.  One difficulty is that,
although these attaching maps are continuous, they are in general not
smooth.  Lemmas~\ref{lemma2.1}, \ref{lemma2.2} and \ref{lemma2.3} show
that the maps possess enough regularity that we can establish a
Smale-like transversality condition
(Definition~\ref{smale_condition}).  Using this definition,
Lemmas~\ref{transversality_lemma} and \ref{finite_trajectories} show
that a finite number of trajectories connect critical points whose
indices differ by one.  Finally a after a careful consideration of
orientations along each trajectory connecting critical points, we
define a $\mbox{degree}$ to each such pair.

\begin{lemma}
Let $K$ be a stratum with dimension $k$ and $L$ be a stratum in the closure of
$K$ with dimension $k-1$.  Let $p\in K$ and suppose that 
$\varphi(s,p)\in K$ for
$s<t_p$ and $\varphi(s,p)\in L$ for $t_p< s\leq T$.  Then there is a 
neighborhood $U\subset K$ of $p$ such that $\varphi(T,\cdot):U\to M$ is 
smooth.
\label{lemma2.1}
\end{lemma}

\begin{proof}
Choose a coordinate system near $\varphi(t_p,p)$ such that for $q\in K$,
$x_{k+1}(q)=\dots=x_n(q)=0$ and for $q\in L$, $x_{k}(q)=\dots=x_n(q)=0$.
We know that the $k$th coordinate satisfies $\varphi_k(t_p,p)=0$, and in 
addition
$$
\left.\frac{\partial \varphi_k}{\partial t}\right|_{(t_p,p)}\neq0,
$$
by the definition of Morse function.

From the Implicit Function Theorem, it follows that there is an open set 
$U\subset K$ containing $p$ and an open set $V\subset \mathbb{R}$ such that 
for each $q\in U$ there is a unique $g(q)\in V$ such that 
$\varphi_k(g(q),q)=0$.  Moreover, the function $g$ is smooth.

We can use this function $g$ to define two additional functions: 
$\gamma:U\to{\mathbb R}\times U$ given by $\gamma(q)=(g(q),q)$ and 
$\tau:\gamma(U)\to{\mathbb R}\times M$ given by $\tau(t,q)=(T-t,\varphi(t,q))$.
Since $\varphi$ is smooth whenever it remains in a single stratum, these are
themselves smooth functions.
We can write 
$$
\varphi(T,q)=(\varphi\circ\tau\circ\gamma)(q),
$$ 
a composition 
of smooth functions, so $\varphi(T,\cdot):U\to M$ is smooth.

\end{proof}

\begin{lemma}
Let $K$ be a stratum with dimension $k-1$ in the closure of a stratum $L$ 
with dimension $k$.  Let $p\in K$ and suppose that 
$\varphi(s,p)\in K$ for
$s\leq t_p$ and $\varphi(s,p)\in L$ for $t_p< s\leq T$.  Then there is a 
neighborhood $U\subset K$ of $p$ such that $\varphi(T,\cdot):U\to M$ is 
smooth.
\label{lemma2.2}
\end{lemma}

\begin{proof}
Choose a coordinate system near $\varphi(t_p,p)$ such that for $q\in K$,
$x_{k}(q)=\dots=x_n(q)=0$ and for $q\in L$, $x_{k+1}(q)=\dots=x_n(q)=0$.
Let $h(s,q)=(-\nabla f)_k(\varphi(s,q))$. Then $h$ is a smooth real valued 
function.  Since the trajectory $\varphi(\cdot,p)$ moves from $K$ to $L$ at
time  $t_p$, we must have $h(t_p,p)=0$.

By the definition of a Morse function, we know that 
$$
\left.\frac{\partial h}{\partial t}\right|_{(t_p,p)}\neq 0.
$$
It then follows from the Implicit Function Theorem, that there is an open set 
$U\subset K$ containing $p$ and an open subset $V\subset \mathbb{R}$ such that 
for each $q\in U$ there is a unique $g(q)\in V$ such that 
$h(g(q),q)=0$.  Moreover, the function $g$ is smooth.

As before, we can define smooth functions 
$\gamma:U\to{\mathbb R}\times U$ given by $\gamma(q)=(g(q),q)$ and 
$\tau:\gamma(U)\to{\mathbb R}\times M$ given by $\tau(t,q)=(T-t,\varphi(t,q))$
and write 
$$
\varphi(T,q)=(\varphi\circ\tau\circ\gamma)(q),
$$
Again, this is a composition of smooth functions, and so
$\varphi(T,\cdot):U\to M$ is smooth.

\end{proof}

Now let $p_2$ be an essential critical point with index $\lambda_2=j$ and 
$p_1$ an essential critical point with index $\lambda_1=j-1$.  Suppose that
$S(p_1)\cap U(p_2)\neq \emptyset$.  Choose a 
$$
q\in S(p_1)\cap U(p_2)\cap M_{c_1+\varepsilon}.
$$
Then there is a $q_1\in U_\varepsilon(p_2)$ and a time $T$ such that 
$q=\varphi(T,q_1)$.

\begin{lemma}
If the trajectory $\varphi(s,q_1)$ only ever changes between strata whose 
difference in dimension is 1, then there is a neighborhood $V$ in the stratum
containing $q_1$ such that $\varphi(T,\cdot):V\to M$ is smooth.
\label{lemma2.3}
\end{lemma}

\begin{proof}
The result follows from repeated application of Lemma~\ref{lemma2.1} and 
Lemma~\ref{lemma2.2}.  Since $f^{-1}(-\infty,r]$ is compact for each choice
of $r$, there can be only 
finitely many stratum changes~(\cite{handron}).  For each stratum 
change, $i=1,\dots,N$, we get a corresponding neighborhood $V_i$.  Since 
there are finitely
many, we may choose
$$
V=V_1\cap\dots\cap V_N.
$$
\end{proof}

\begin{definition}
We say that the Morse function $f$ is a {\em Morse-Smale function} if for each 
pair of critical points, $p$ and $p'$ with $f(p)>f(p')$, and for 
$\varepsilon$ sufficiently small
the parametrized space $U(p)$ and the set 
$S_\varepsilon(p')=S(p')\cap M_{f(p')+\varepsilon}$ 
intersect transversely in the following sense: if $q\in
S_\varepsilon(p')\cap U(p)$, $K$ is the stratum containing $q$, and $q$ has a
neighborhood $N\subset M$ such that $S_\varepsilon(p')\cap U(p)\cap
N=S_\varepsilon(p')\cap U(p)\cap K$, then $U(p)$ and
$S_\varepsilon(p')\cap K$ intersect transversely in $K$.
\label{smale_condition}
\end{definition}

Let $p_i$ and $p_j$ be two critical points with indices $\lambda_i$ and 
$\lambda_j$, respectively, satisfying  $\lambda_j-\lambda_i=1$.  Let
$\mathfrak{T}$ be the set of trajectories from $p_j$ to $p_i$.  We wish
to show that the
transversality condition in Definition~\ref{smale_condition} together 
with the compactness of $M$ ensure that $\mathfrak{T}$
is a finite set.

We begin with the following

\begin{lemma}
\label{transversality_lemma}
Suppose that $\widetilde{M}$ is an $m$-manifold, $\widetilde{S}$ a $k$-dimensional regular 
submanifold.  Let $\widetilde{U}$ be an $(m-k)$-manifold, and 
$\Phi:\widetilde{U}\to\widetilde{M}$ a smooth map such that $\widetilde{S}\pitchfork\Phi(\widetilde{U})$.  Then 
$\Phi^{-1}(\widetilde{S}\cap\Phi(\widetilde{U}))$ consists of isolated points in $\widetilde{U}$.
\end{lemma}

\begin{proof}
First, since $\dim(\widetilde{S})+\dim(\Phi(\widetilde{U}))=k+(m-k)=m$, the transverse intersection
of $\widetilde{S}$ and $\Phi(\widetilde{U})$ is zero dimensional.  Let $q\in \widetilde{U}$ such that 
$\Phi(q)\in \widetilde{S}$.  We will show that there is a 
neighborhood $N_q$ of $q$ such that $\Phi^{-1}(\widetilde{S})\cap N_q=\{q\}$.

let $\mathbf{x}=(x^1,\dots,x^m)$ be a coordinate chart on $\widetilde{M}$ 
such that 
$\mathbf{x}(\Phi(q))=(0,\dots,0)$, and $(x^1,\dots,x^k)$ gives a coordinate 
chart on $\widetilde{S}$.  Then for $p\in \widetilde{S}$
$$
x^i(p)=0,\quad\mbox{for }i=k+1,\dots,m.
$$

Let $\bar{\mathbf{x}}=(x^{k+1},\dots,x^m)$.  We want to show that 
$\tilde{\mathbf{x}}=\bar{\mathbf{x}}\circ\Phi$ is a coordinate chart on $\widetilde{U}$.
It suffices to show that if $\mathbf{y}$ is a coordinate chart on $\widetilde{U}$, then
$D(\tilde{\mathbf{x}}\circ\mathbf{y}^{-1})$ is non-singular.

Since $\widetilde{S}\pitchfork\Phi(\widetilde{U})$, the vectors
$$
\left\{
\frac{\partial}{\partial x^1}|_{\Phi(q)},
\dots,
\frac{\partial}{\partial x^k}|_{\Phi(q)},
\Phi_\ast\left(\frac{\partial}{\partial y^{k+1}|_q}\right),
\dots,
\Phi_\ast\left(\frac{\partial}{\partial y^{n}|_q}\right)
\right\}
$$
form a basis for $T_{\Phi(q)}(\widetilde{M})$, as do 
$\left\{
\frac{\partial}{\partial x^1}|_{\Phi(q)},
\dots,
\frac{\partial}{\partial x^m}|_{\Phi(q)}
\right\}$.
We can thus write
$$
\Phi_\ast\left(\frac{\partial}{\partial y^{i}}\right)=\sum_{j=1}^n
\alpha_i^j\frac{\partial}{\partial x^j}\quad\mbox{ for }i=k+1,\dots,n.
$$

Using this expression, it is clear that the matrix
$$
\begin{bmatrix}
1              & 0 & \cdots & 0              \\
0              & 1 & \cdots & 0              \\
\vdots         &   &        & \vdots         \\
\alpha_{k+1}^1 &   & \cdots & \alpha_{k+1}^m \\
\vdots         &   &        & \vdots         \\
\alpha_{m}^1 &   & \cdots & \alpha_{m}^m 
\end{bmatrix}
=
\left[
\begin{array}{c|c}
I_k  & 0 \\ \hline
\ast & 
\begin{matrix}
\alpha_{k+1}^{k+1} & \cdots & \alpha_{k+1}^m \\
\vdots             & \ddots & \vdots         \\
\alpha_m^{k+1}     & \cdots & \alpha_m^m     
\end{matrix}
\end{array}
\right]
=
\left[
\begin{array}{c|c}
I_k  & 0 \\ \hline
\ast & M
\end{array}
\right]
$$
is non-singular (having linearly independent rows).  It follows that the 
matrix $M$ is non-singular.  But $M$ is the transpose of the matrix for
$$
D(\bar{\mathbf{x}}\circ(\Phi\circ\mathbf{y}^{-1}))
=D((\bar{\mathbf{x}}\circ\Phi)\circ\mathbf{y}^{-1})
=D(\tilde{\mathbf{x}}\circ\mathbf{y}^{-1}),
$$
evaluated at $\mathbf{y}(q)$, so 
$D(\tilde{\mathbf{x}}\circ\mathbf{y}^{-1})$ is non-singular in a neighborhood
of $\mathbf{y}(q)$.

Now, since $\tilde{\mathbf{x}}$ is a coordinate chart in a neighborhood of 
$q$, there is an open set
$N_q$ of $q$ such that for $q'\in N_q$ with $q'\neq q$, 
$\tilde{\mathbf{x}}(q')\neq(0,\dots,0)$.  So $\Phi(q')\notin S$.

\end{proof}

\begin{lemma}
Let $f:M\to\mathbb{R}$  be a Morse-Smale function on a manifold with corners 
$M$.
Let $p_i$ and $p_j$ be two critical points of $f$ with indices $\lambda_i$ and 
$\lambda_j$ respectively, satisfying $\lambda_j-\lambda_i=1$.  
Then the set
$\mathfrak{T}$ of trajectories from $p_j$ to $p_i$ is finite.
\label{finite_trajectories}
\end{lemma}

\begin{proof}
Let $c_i=f(p_i)$ and $c_j=f(p_j)$.  Then $\partial U_\varepsilon(p_j)$
is diffeomorphic to $S^{\lambda_j-1}$.  For each element of
$\tau\in\mathfrak{T}$ there is exactly one $q_\tau\in \partial
U_\varepsilon(p_j)\cap\tau$.

Suppose that $\mathfrak{T}$ is an infinite set.  Then, since $\partial U_\varepsilon(p_j)$ is
compact, there is a convergent sequence $\{q_i\}_{i=1}^\infty$ such that 
$q_i\in \partial U_\varepsilon(p_j)\cap\tau_i$ for some $\tau_i\in\mathfrak{T}$.  Let
$q_0=\lim_{i\to\infty}q_i$.

We wish to show that $q_0=\partial U_\varepsilon(p_j)\cap\tau_0$ for some 
$\tau_0\in\mathfrak{T}$.  Choose $T$ such that for $p= \varphi(T,q_0)$,
$f(p)<c_i+\varepsilon$.  Then for $i$ sufficiently large
$f(\varphi(T,q_i))<c_i+\varepsilon$, and since for each $i$
$\lim_{t\to\infty}\varphi(t,q_i)=p_i$, $\varphi(T,q_i)\in S_\varepsilon(p_i)$
for $i$ sufficiently large.  Now, $\varphi$ is continuous, and 
$S_\varepsilon(p_i)$ is closed, so 
$p=\lim_{i\to\infty}\varphi(T,q_i)\in S_\varepsilon(p_i)$.
It follows that $q_0$ lies on some trajectory connecting $p_j$ to $p_i$.

Now we will show that $p_0$ is isolated in $\partial
U_\varepsilon\cap\mathfrak{T}$, contradicting the choice of $p_0$.
Let $K$ be the stratum containing $p=\varphi(T,q_0)$ and set $c=f(p)$.
Let $\widetilde{M}=f^{-1}(c)\cap K$,
$\widetilde{S}=S_\varepsilon(p_i)\cap\widetilde{M}$.  Assume there is
a neighborhood $U_0$ of $q_0$ in $M$ on which $\varphi(T,\cdot)$ is
smooth.  (If there is not, replace $T$ with a slightly larger $T'$.)
For each $q\in U_0$ there is a trajectory of $\varphi$ containing $q$.
We will denote this trajectory by $\gamma_q$.  For a given $q\in U_0$,
let $\tilde{q}=\gamma_q\cap\partial U_\varepsilon(p_j)$,
$\bar{q}=\gamma_q\cap\varphi(T,\cdot)^{-1}(\widetilde{S})$, and let
$t_q\in\mathbb{R}$ such that
$$
\begin{array}{c}
\varphi(t_q,\tilde{q})=\bar{q}\quad\mbox{ if }f(\tilde{q})\geq f(\bar{q})  \\
\varphi(-t_q,\bar{q})=\tilde{q}\quad\mbox{ if }f(\tilde{q})\leq f(\bar{q}).
\end{array}
$$

The maps $q\mapsto\tilde{q}$, $q\mapsto\bar{q}$ and $q\mapsto t_q$ are
all smooth, and so the function $\Phi_0:U_0\to M$ given by
$\Phi_0(q)=\varphi(T+t_q,q)$ is smooth.  Let $\widetilde{U}=\partial
U_\varepsilon(p_j)\cap U_0$, and $\Phi=\Phi_0\vert_{\widetilde{U}}:\widetilde{U}\to\widetilde{M}$.

We can let $k=\dim(\widetilde{S})=\dim(K)-\lambda_i-1$ and 
$m=\dim(\widetilde{M})=\dim(K)-1$.  Then $\dim(\partial U_\varepsilon(p_j))=\lambda_j-1=m-k$.
Choose a coordinate chart at $p\in\widetilde{M}$, $\mathbf{x}=(x^1,\dots,x^m)$,
such that $(x^1,\dots,x^k)$ is a coordinate chart on $\widetilde{S}$ and
$$
x^i(p')\neq0\quad\mbox{ whenever }p'\notin\widetilde{S}\mbox{ and }i>k.
$$
This can be done in such a way that 
$\frac{\partial}{\partial x^{k}}\vert_p$ is tangent to the trajectory
$\tau_0$.

Now choose a coordinate chart $\mathbf{y}=(y^{k+1},\dots,y^{m+1})$  
on $U_0$, such that $\mathbf{y}=(y^{k+1},\dots,y^m)$ is a 
coordinate chart on $\widetilde{U}$ and 
$\frac{\partial}{\partial y^{m+1}}\vert_{q_0}$ is tangent to the trajectory 
$\tau_0$.  Then 
$$
\Phi_\ast\left(\frac{\partial}{\partial y^i}\vert_{q_0}\right)={\Phi_0}_\ast\left(\frac{\partial}{\partial y^i}\vert_{q_0}\right)\in T_p\widetilde{M}
\mbox{ for }i=k,\dots,m.
$$
In addition, ${\Phi_0}_\ast(\frac{\partial}{\partial y^n}\vert_{q_0})$ will be
tangent (at $p$) to the trajectory $\tau_0$.

Since $f$ is a Morse-Smale function, the vectors
$$
\left\{
\frac{\partial}{\partial x^{1}}\vert_p,
\dots,
\frac{\partial}{\partial x^{k+1}}\vert_p,
{\Phi_0}_\ast\left(\frac{\partial}{\partial y^{k+1}}\vert_{q_0}\right),
\dots,
{\Phi_0}_\ast\left(\frac{\partial}{\partial y^{m+1}}\vert_{q_0}\right)
\right\}
$$
span $T_pK$.  Since the vectors $\frac{\partial}{\partial x^{k+1}}\vert_p$
and ${\Phi_0}_\ast(\frac{\partial}{\partial y^{m+1}}\vert_{q_0})$ are tangent 
to $T_p\tau_0$, and the remaining $m=\dim(\widetilde{M})$ vectors are all in 
$T_p\widetilde{M}$,
$$
\left\{
\frac{\partial}{\partial x^{1}}\vert_p,
\dots,
\frac{\partial}{\partial x^{k}}\vert_p,
\Phi_\ast\left(\frac{\partial}{\partial y^{k+1}}\vert_{q_0}\right),
\dots,
\Phi_\ast\left(\frac{\partial}{\partial y^{m}}\vert_{q_0}\right)
\right\}
$$
form a basis for $T_p\widetilde{M}$.  The first $k$ vectors of this basis
are a basis of $T_p\widetilde{S}$, and the remaining $m-k$ vectors are a 
basis for $T_p\Phi(\partial U_\varepsilon(p_j))$.  It follows that
$$
\widetilde{S}\pitchfork\Phi(\partial U_\varepsilon(p_j)).
$$

By Lemma~\ref{transversality_lemma}, 
$\Phi^{-1}(\widetilde{S}\cap\Phi(\widetilde{U}))$ consists of isolated points.
Note that for sufficiently large $i$, 
$q_i\in\Phi^{-1}(\widetilde{S}\cap\Phi(\widetilde{U}))$ and in addition 
$q_0\in\Phi^{-1}(\widetilde{S}\cap\Phi(\widetilde{U}))$.  This, however,
contradicts the fact that $\lim_{i\to\infty}q_i=q_0$.  Consequently, the 
number of trajectories connecting $p_j$ to $p_i$ must be finite.
\end{proof}

Now assume that $M$ is an oriented manifold.  
At each essential critical point $p$, the
tangent space $T_pM$ can be decomposed into a stable space and an unstable 
space, $T_pM=E^+(p)\oplus E^-(p)$.  Choose an orientation for the 
unstable space $E^-(p)$.  This induces an orientation on $E^+(p)$ in the 
following way:

Let $\lambda_p$ be the index of $p$, and choose a basis
$v_1,\dots,v_{\lambda_p}\in E^-(p)$ which represents the orientation
on $E^-(p)$.  Then choose $v_{\lambda_p+1},\dots,v_n\in E^+(p)$ such
that $v_1,\dots,v_n$ is a basis which represents the orientation on
$T_pM$.  The vectors $v_{\lambda_p+1},\dots,v_n$ determine an
orientation for $E^+(p)$.

Now for critical points $p_i$ and $p_j$, the orientation on $E^+(p_i)$ 
determines an orientation on 
$S_\varepsilon(p_i)$, which 
in turn determines an orientation on $T_qS(p_i)$ for $q\in S(p_i)$.  In a 
similar way, the 
choice of orientation on $E^-(p_j)$ induces an orientation on 
$T_{q'}U(p_j)$ for $q'\in U_\varepsilon(p_j)$.

Note that if the function $f:M\to{\mathbb R}$ is a Morse-Smale function
and $\varphi(T,q')=q$, then
the derivative $\varphi(T,\cdot)_\ast$ restricted to $U(p_j)$ has full rank at
$\varphi(T,\cdot)^{-1}(q)$.

Now the map $\varphi(T,\cdot)_\ast$ pushes forward the orientation on 
$T_{q'}U(p_j)$ to give an
orientation on $T_q\varphi(T,U(p_j))$.  Since 
$\dim(\varphi(T,U(p_j)))=j$ and 
$\dim(S_\varepsilon(p_i))=n-j+1$, these two sets intersect in the
one dimensional trajectory through $q$.  We can choose a basis for 
$T_q\varphi(T,U(p_j))$ which represents the induced orientation and such that
the last vector is $\frac{d}{dt}|_T\varphi(t,q_1)$.  We also may choose 
a basis representing the orientation on $T_qS(p_i)$, with 
$\frac{d}{dt}|_T\varphi(t,q_1)$ as its first 
vector.  When combined, these bases give an orientation for $T_qM$.  
If this orientation agrees with the chosen orientation on $M$ then we say
the sign of the trajectory from $p_j$ to $p_i$ through $q_1$ is $+1$.  
Otherwise, the sign is $-1$.

Now if $\mathfrak{T}$ is the set of trajectories from $p_j$ to $p_i$, 
Lemma~\ref{finite_trajectories} ensures that $\mathfrak{T}$
is a finite set.  So, we can define
$$
\mbox{degree}(p_j,p_i)=\sum_{\tau\in\mathfrak{T}}\mbox{sign}(\tau).
$$
This degree, we shall see, is the degree of the attaching map 
$\Sigma_j$ at $p_i$.

\section{The Morse Complex for a Morse Function on a Manifold with Corners}
\label{complex}

\subsection{The Geometric Complex of a Morse Function}

Let $P_k$ denote the cell $U_\varepsilon(p_k)$.
We know that $M$ is homotopy equivalent to the cellular complex 
$$
X=P_1\cup_{\Sigma_2}P_2\cup\dots\cup_{\Sigma_m}P_m.
$$

Recall that $\Sigma_j$ was defined by
$$
\Sigma_j=\Psi_2(1,\cdot)\circ\Phi_2(1,\cdot)\circ\dots\circ\Psi_j(1,\cdot)\circ\Phi_j(1,\cdot)\circ\sigma_j,
$$
which attaches $P_j$ to
$P_1\cup_{\Sigma_2}P_2\cup\dots\cup_{\Sigma_{j-1}}P_{j-1}$.  The map
$\Sigma_{j}$ is defined by a sequence of maps derived from flows.
They alternate between a $\Phi_k$ which preserves the
$\varphi$-trajectories, and a $\Psi_k$ which does not.

We wish to show that this cellular complex is, in fact, a CW-complex.
In order to do so we need the following result regarding the behavior
of the trajectories of $\varphi$.

\begin{lemma}
Let $p_i$ and $p_j$ be essential critical points of a Morse-Smale
function $f$ on a manifold with corners $M$.  If $p_j\in\overline{U(p_i)}$
then $\mbox{index}(p_j)<\mbox{index}(p_i)$.  If $p_i\in\overline{S(p_j)}$
then $\mbox{index}(p_i)>\mbox{index}(p_j)$. 
\label{critical_point_lemma}
\end{lemma}

This lemma allows the maps $\Phi_k$ and $\Psi_k$ to be defined in a
way that preserves $\varphi$-trajectories connecting critical points
whose indices differ by one.

\begin{proof}
We will show that the boundaries of $U(p_i)$ and of $S(p_j)$ are
composed of a union of trajectories.  Suppose, in this case, that
$p_j\in\overline{U(p_i)}$.  Then there is a trajectory $\tau\in
S(p_j)$, with $\tau\subset\overline{U(p_i)}$.  The upper endpoint of
$\tau$ must also be a critical point $p_k$ in $\overline{U(p_i)}$.
Since $f$ is a Morse-Smale function, it follows that
$\mbox{index}(p_k)>\mbox{index}(p_j)$. 

Traveling up from $p_j$ in this
way, there is a sequence of essential critical points with increasing
indices.  Since the topmost critical point in this sequence must be
$p_i$, the result follows.  The second half of the lemma is
established in a similar way.

To show that the boundary of $U(p_i)$ and of $S(p_j)$ is composed of a
union of trajectories, we will assume that $q\in\overline{U(p_i)}$ and
we will show two things.  First, that for any $t_0>0$,
$\varphi(t_0,q)\in\overline{U(p_i)}$ and second, that that there exist
$t_0$ and $q'$ such that $q'\in\overline{U(p_i)}$ and
$\varphi(t_0,q')=q$.

If $q\in\overline{U(p_i)}$, then there is a sequence $\{q_n\}\subset
U(p_i)$ such that $\lim_{n\to\infty}q_n=q$.  Since $\varphi(t_0,\cdot)$ is
continuous, it follows that
$\lim_{n\to\infty}\varphi(t_0,q_n)=\varphi(t_0,q)$.

The second assertion is more of a challenge, since $\varphi$ provides
a flow in only one direction.  Again, let $\{q_n\}\subset
U(p_i)$ be a sequence such that $\lim_{n\to\infty}q_n=q$. Let $r=f(q)$.

If $q$ is an essential critical point, then let $q'=q$ and $t_0$ be
any real number greater than zero.  Then $\varphi(t_0,q')=q$.  

If $q$ is not an essential critical point, let then for each $n$, let
$q_n'$ denote the intersection of the trajectory containing $q_n$ with
$f^{-1}(r)$.  Then $\lim_{n\to\infty}q_n'=q$ also.

Let $r'>r=f(q)$ be such that $f^{-1}([r,r'])$ contains no critical
points (at least in a neighborhood of $q$).  Then let $q_n''$ denote
the intersection of the trajectory containing $q_n$ with $f^{-1}(r')$.
Since $f^{-1}(r')$ is a closed submanifold, it is compact.  There is a
subsequence of $\{q_n''\}$ that converges to a point $q'$.  Relabel,
if necessary, so that $\lim_{n\to\infty}q_n''=q'$.  Note that
$q'\in\overline{U(p_i)}$.

The map from $f^{-1}(r')$ to $f^{-1}(r)$ given by allowing each point
of $f^{-1}(r')$ to flow by $\varphi$ until it hits $f^{-1}(r)$ is both
well defined and continuous (at least in a neighborhood of
$q'$). Since this map carries each $q_n''$ to $q_n'$, it follows that
the limit $q'$ is carried to the corresponding limit $q$.  Thus there
is some $t_0$ such that $\varphi(t_0,q')=q$.

It follows that $\overline{U(p_i)}$ contains the entirety of some
trajectory passing through $q$.

The above argument uses only the fact that $\overline{U(p_i)}$ is a
union of trajectories.  The argument applies equally well to
$\overline{S(p_j)}$

\end{proof}

Note that this lemma does not imply that if $\tau$ is a trajectory
that intersects $\overline{U(p_i)}$, then
$\tau\subset\overline{U(p_i)}$.  Because trajectories combine wherever
$-\nabla f$ points outward, a point may lie in more than one
trajectory.  Of course, all trajectories containing a particular point
must have the same lower limit.

\begin{lemma}
The cellular complex $X$ of the Morse-Smale function $f:M\to\mathbb{R}$
is a CW-complex.
\label{cw-complex_lemma}
\end{lemma}

\begin{proof}
By construction, $\Sigma_i$ maps $\partial P_i=\partial
U_\varepsilon(p_i)$ onto a union of cells.  The important point here
is to show that these cells have dimension smaller than that of
$P_i$.  Ideally, a cell $P_j$ will be in the image of $\Sigma_i$ if
and only if there is a $\varphi$-trajectory connecting $p_i$ to $p_j$.

By Lemma~\ref{critical_point_lemma}, for each critical point $p_i$
there is an open set $W$ containing $p_i$, such that $U(p_j)\cap
W=\varnothing$ whenever $\mbox{index}(p_j)\leq\mbox{index}(p_i)$ and
$S(p_j)\cap W=\varnothing$ whenever
$\mbox{index}(p_j)\geq\mbox{index}(p_i)$.  We can then choose
$\varepsilon$ small enough that $\Psi_{i+1}$ is the identity outside
$W$.

Now suppose that $S(p_j)\cap U(p_i)=\varnothing$, but
$P_j\subset\mbox{image}(\Sigma_i)$.  The sets $S(p_j)$ and $U(p_i)$
must have been ``confused'' in a neighborhood of some essential
critical point $p_k$ by $\Psi_{k+1}$.  It follows that
$\mbox{index}(p_k)<\mbox{index}(p_i)$ and
$\mbox{index}(p_j)<\mbox{index}(p_k)$. It follows that
$$
\mbox{index}(p_j)\leq\mbox{index}(p_i)-2$$.

Hence, the boundary of each cell is mapped to a union of lower
dimensional cells, and the geometric complex is a CW-complex.

\end{proof}

The proof gives us some other information about the attaching maps of
this CW-complex.  Suppose that $\Sigma$ is the attaching map for
$P=U_\varepsilon(p)$ and $Q=U_\varepsilon(q)$ is a cell with dimension
one less than that of $P$.  Then $\Sigma^{-1}(q)$ is exactly the set
of points in $\partial U_\varepsilon(q)$ that lie on trajectories
connecting $p$ to $q$.

This CW-complex gives rise to a chain complex
$$
0 \longleftarrow  W_0  \overset{\partial_1}{\longleftarrow}  W_1  \overset{\partial_2}{\longleftarrow} \dots \overset{\partial_n}{\longleftarrow}  W_n.
$$
The free Abelian groups $W_i$ are generated by the $i$-cells of the
above CW-complex.  We will refer to this chain complex as the {\em
geometric complex} of $f:M\to\mathbb{R}$.

\subsection{The Algebraic Complex of a Morse Function}
Let $C_j(f)$ be the set of essential critical points of $f$ with index
$j$, and $V_j$ the free Abelian group of formal ${\mathbb Z}$-linear
combinations of the elements of $C_j(f)$.  For each
$j=1,\dots,n=\dim(M)$ define a map $\partial'_j:V_j\to V_{j-1}$ by
setting
$$
\partial'_j(p)=\sum_{q\in C_{j-1}(f)}\mbox{degree}(p,q)q
$$
for each $p\in C_j(f)$, and extending linearly to all of $V_j$.  For
$j=0$, set $\partial'_0(p)=0$

In the next section, we will see that the free Abelian groups $V_j$
and maps $\partial'_j$ form a chain complex
$$
0 \longleftarrow  V_0  \overset{\partial'_1}{\longleftarrow}  V_1  \overset{\partial'_2}{\longleftarrow} \dots \overset{\partial'_n}{\longleftarrow}  V_n
$$
We will refer to this chain complex as the {\em algebraic complex}
of $f:M\to\mathbb{R}$.

\subsection{Main Theorem}

We are now ready to state and prove the main result of this article.

\begin{theorem}
Let $M$ be a compact orientable $n$-manifold with corners, and 
$f:M\to\mathbb{R}$ be a Morse-Smale function satisfying the conditions of 
Lemma~\ref{lemma2.3}.  Then the algebraic complex of
$f:M\to\mathbb{R}$ is a chain complex.  Moreover, the algebraic
complex and geometric complex have the same homology groups. 
\label{main_theorem}
\end{theorem}

\begin{proof}
Let $C_k(X)$ be the set of cells in $X$ having dimension $k$.  As
above, $C_k(f)$ is the set of critical points of $f$ having index
$k$.  These sets are in one to one correspondence;  each $P\in C_k(X)$
satisfies $P=U_\varepsilon(p)$ for some $p\in C_k(f)$ and for each
$p\in C_k(f)$, $P=U_\varepsilon(p)\in C_k(X)$.

The bijection $h_k:C_k(X)\to C_k(f)$ by $U_\varepsilon(p)\mapsto p$
can be extended to an isomorphism $h_k:W_k\to V_k$ by linearity.
Thus we have a chain map $h$ given by
$$
\begin{CD}
0    @<<<  W_0   @<{\partial_1}<<   W_1   @<{\partial_2}<<   \dots   @<{\partial_n}<<  W_n   @<<< 0\\
@VVV       @V{h_0}VV                @V{h_1}VV                @VVV                        @V{h_n}VV \\
0    @<<<  V_0   @<{\partial'_1}<<  V_1   @<{\partial'_2}<<  \dots   @<{\partial'_n}<<  V_n  @<<< 0
\end{CD}
$$
Where each $h_i$ is an isomorphism.  It suffices to show that
$\partial_k=h_{k-1}^{-1}\circ\partial'_k\circ h_k$.

The map $\partial_k$ is defined by
$$
\partial_k(P_j)=\sum_{P_i\in C_{k-1}(X)}\deg(\Sigma_j,P_i)P_i
$$
for each $P_j\in C_k(X)$, where $\Sigma_j$ is the attaching map for the 
cell $P_j$.  
The degree $\deg(\Sigma_j,P_i)$ is determined by choosing a $q\in P_i$ such that 
$\Sigma_j(\partial P)$ intersects $P_i$ transversely at $q$, and then adding the
number of points in $\Sigma_j^{-1}(q)$ at which $\Sigma_j$ is orientation
preserving, and subtracting the number of points in $\Sigma_j^{-1}(q)$ at 
which $\Sigma_j$ is orientation reversing.  Given $P_j\in C_k(X)$ and 
$P_i\in C_{k-1}(X)$ we must determine this degree.

We have the relationship
$$
\deg(\Sigma_j,P_i)=\deg(\Sigma_{ij},U_\varepsilon(p_i)).
$$

Consider a point $q\in\partial U_\varepsilon(p_j)\cap\tau$  for some
trajectory connecting $p_j$ to $p_i$.  The map $\Sigma_{ij}$ is orientation 
preserving at $q$ if $\mbox{sign}(\tau)=+1$, and orientation reversing if 
$\mbox{sign}(\tau)=-1$.  It follows that 
$$
\deg(\Sigma_j,P_i)=\deg(\Sigma_{ij},p_i)=\mbox{degree}(p_j,p_i).
$$
and so
$$
\deg(\Sigma_j,P_i)=\mbox{degree}(h_k(p_j),h_{k-1}(p_i)).
$$
It follows that 
$$
\partial_k(P)=\sum_{P_i\in
C_{k-1}(X)}\mbox{degree}(h_k(p_j),h_{k-1}(p_i))h_{k-1}^{-1}(p_i). 
$$
and so $\partial_k=h_{k-1}^{-1}\circ\partial'_k\circ h_k$ as claimed.
This proves the theorem.

\end{proof}

It should be pointed out that in the the descending set of an
essential critical point may not be homeomorphic to an open disk.
Consequently, the descending sets may not be interpreted directly as
cells in a CW-complex, as they may in the classical theory.
Nonetheless, this theorem provides a useful parallel with the classical
Morse theory in the setting of manifolds with corners.

%\section{Acknowledgments}
%I would like to that my editor, reviewers and all those whose advice 
%assisted in writing and publishing this article.

\end{document}